\documentclass[11pt]{amsart}
\usepackage{amsfonts}
\usepackage{mathrsfs}
\usepackage{graphicx}
\usepackage{booktabs}
\usepackage{array}
\usepackage{amsmath}\allowdisplaybreaks[4]
\usepackage [latin1]{inputenc}
\usepackage[numbers,sort&compress]{natbib}

\newtheorem{thm}{Theorem}[section]

\newtheorem{lem}[thm]{Lemma}
\newtheorem{prop}[thm]{Proposition}
\theoremstyle{definition}

\theoremstyle{remark}
\newtheorem{rem}[thm]{Remark}

\numberwithin{equation}{section}

\begin{document}

\title[On the global solvability of the generalised Navier-Stokes system]{{On the global solvability of the generalised Navier-Stokes system in critical Besov spaces }}
\author[Huiyang Zhang, Shiwei Cao \& Qinghua Zhang]
{Huiyang Zhang, Shiwei Cao \& Qinghua Zhang*} 

\thanks{{\it Address of the Corresponding Author:} School of Mathematics and Statistics,
  Nantong University,
6 Seyuan Road, Nantong City 226019, Jiangsu Province, P. R. China.}
\thanks{ E-mail: zhangqh@ntu.edu.cn}


\begin{abstract}
This paper is devoted to the global solvability of the Navier-Stokes system with fractional Laplacian $(-\Delta)^{\alpha}$ in $\mathbb{R}^{n}$ for $n\geq2$, where the convective term has the form $(|u|^{m-1}u)\cdot\nabla u$ for $m\geq1$. By establishing the estimates for the difference $|u_{1}|^{m-1}u_{1}-|u_{2}|^{m-1}u_{2}$ in homogeneous Besov spaces, and employing the maximal regularity property of $(-\Delta)^{\alpha}$ in Lorentz spaces, we prove global existence and uniqueness of the strong solution of the Navier-Stokes  in critical Besov spaces for both $m=1$ and $m>1$.
\end{abstract}
\maketitle
\section{Introduction and main results}
This paper deals with the generalised Navier-Stokes system (or symbolically (GNS)) in $\mathbb{R}^{n}$ ($n\geq2$), namely
\begin{equation}\label{eqn:gns}
\left\{\begin{array}{l}
\partial_{t}u+\mu(-\Delta)^{\alpha}u+J_{m}(u)\cdot\nabla u+\nabla \pi=f(x,t),\;\;t>0,\;x\in\mathbb{R}^{n};\\
\textrm{div}u=0,\;\;t>0,\;x\in\mathbb{R}^{n};\\
\;u(0,x)=a(x),\;\;x\in\mathbb{R}^{n},
\end{array} \right.
\end{equation}
where $(-\Delta)^{\alpha}$ is the fractional power of the negative Laplacian for $\alpha>0$ and $J_{m}(u)=|u|^{m-1}u$ for $m\geq1$. We use $u(x,t)=(u_{1}(x,t),u_{2}(x,t),\cdots,u_{n}(x,t))$ and $\pi(x,t)$ to denote the velocity and inner pressure respectively of the flow at the place $x\in \mathbb{R}^{n}$ and the time $t>0$.  While $a(x)=(a_{1}(x),a_{2}(x),\cdots,a_{n}(x))$ represents the initial velocity of the flow. For the sake of simplicity, here the viscosity coefficient $\mu$ is commonly assumed to $1$.

Evidently, if $\alpha=m=1$, then (\ref{eqn:gns}) reduces to the standard incompressible Navier-Stokes system, or (NS) in symbol. The classic result is that if $a\in L^{2}$ (here notation $(\mathbb{R}^{n})$ after $L^{2}$ is omitted for convenience), then (NS) has a globally existing Leray-Hopf weak solution (refer to \cite{ho1951,ms1988,ga2000}), while uniqueness and smoothness of the weak solution for $n\geq3$ is a well-known open problem.  Recall that if the initial velocity has more regularity, then the strong solution of (NS) exists locally. Moreover, if $a$ lies in a critical spaces such as $\dot{H}^{n/2-1}$ and $L^{n}$ etc with sufficiently small size, then (NS) permits a globally existing strong solution (see \cite{fk1964,ka1984,ko1989,iw1989} for references).

If $m=1$ merely, then for $\alpha\geq(n+2)/4$, (GNS) has a unique global regular solution. This well-known result was revealed  by Lions \cite{li1969}.
As for $\alpha<(n+2)/4$, well-posedness for (GNS) reserves only for small data (small norm or small lifespan). By the analytic semigroup approach in Besov spaces, Chen \cite{ch2017} established local existence of the Gevrey regular solution in the class $L^{\infty}(0,T;\dot{B}^{1-2\alpha+n/p}_{p,\infty})$ for $\alpha>1/2$ and $1<p<\infty$ with $2-2\alpha+n/p>0$, and strong solution in the class $C([0,T];\dot{B}^{n/p}_{p,\infty})$ for $\alpha=1/2$ and $1<p<\infty$. Duan \cite{du2018} addressed the global existence of the strong solution to $3D$ (GNS) in homogeneous Sobolev space with small initial data in $\dot{H}^{(5-4\alpha)/2}$ for $1/2<\alpha<5/4$.  In Liu \cite{li2020}, by laying the small assumption on the norm of $e^{-t(-\Delta)\alpha}a\cdot\nabla e^{-t(-\Delta)\alpha}a$ in $L^{1}(0,\infty;\dot{B}^{1-2\alpha+n/p}_{p,1})$ instead of the small assumption on $\|a\|_{\dot{B}^{1-2\alpha+n/p}_{p,1}}$, global solvability of (GNS) in the homogeneous Besov space $\dot{B}_{\infty,\infty}^{-(2\beta-1)}$ for $1/2<\alpha\leq1$ was established.
Other discussions on the same theme can also be founded in \cite{wu2005,yz2014} etc.

This paper mainly deals with the case $\alpha>1/2$ and $m>1$. Since here the generalised  convective term $J_{m}(u)\cdot\nabla u$  could not be rewritten to $\textrm{div}(J_{m}(u)\otimes u)$ anymore, equality $\langle J_{m}(u)\cdot\nabla v,J_{p}(v)\rangle=0$ does not hold for any smooth function $v$ and $p\geq1$. So traditional method of spectral localization performed on (GNS) to get the a prior estimates for the strong solution seems not feasible. For this reason, we given up this method, and turn to appeal to the semigroup $e^{-t\mu(-\Delta)^{\alpha}}$. With some revisions to the results established in \cite{ks2019}, Lorentz maximal regularity of $e^{-t\mu(-\Delta)^{\alpha}}$ on homogeneous Besov spaces will be proved firstly. After that, local H\"{o}lder continuity of the operator $J_{m}(f)$ and estimates for $J_{m}(u_{1})\cdot\nabla u_{1}-J_{m}(u_{2})\cdot\nabla u_{2}$ will be investigated into two cases: $1<m<2$ and $m\geq2$.
Using them as  bases, global existence of the strong solution to (GNS) is finally proved. To our best knowledge, this is the first attemptation to deal with the generalised Navier-Stokes equations with the convective term $J_{m}(u)\cdot\nabla u$ in Besov spaces.

Before stating the main results, we first make a brief review on the critical function spaces. It is easy to check that, system \eqref{eqn:gns} is invariant under the following scaling transformations
\begin{align}\label{tra:sca}
&u_{\lambda}(x,t)=\lambda^{(2\alpha-1)/m}u(\lambda x,\lambda^{2\alpha}t),\;\pi_{\lambda}(x,t)=\lambda^{(2\alpha-1)/m+2\alpha-1}\pi(\lambda x,\lambda^{2\alpha}t)\;\mathrm{and}\\
 \nonumber  &f_{\lambda}(x,t)=\lambda^{(2\alpha-1)/m+2\alpha}f(\lambda x,\lambda^{2\alpha}t).
\end{align}
Since the norm $\|u(\lambda\cdot)\|_{\dot{B}_{p,r}^{s}}$ is equivalent to $\lambda^{s-n/p}\|u\|_{\dot{B}_{p,r}^{s}}$ for all $\lambda>0$, one can easily check that $\|a\|_{\dot{B}_{p,r}^{n/p-(2\alpha-1)/m}}$ is invariant under the transformation $a\mapsto\lambda^{(2\alpha-1)/m}a(\lambda \cdot)$, while $\|u\|_{L^{\rho,r}(0,\infty;\dot{B}_{p,1}^{n/p-(2\alpha-1)/m+2\alpha/\rho})}$ is invariant under the transformation $u\mapsto u_{\lambda}$. In this sense, $\dot{B}_{p,r}^{n/p-(2\alpha-1)/m}$ and  $L^{\rho,r}(0,\infty;\dot{B}_{p,1}^{n/p-(2\alpha-1)/m+2\alpha/\rho})$ are called critical initial space and critical temporal-spatial space respectively.

For the sake of convenience, throughout this paper, notations $L^{\rho,r}(0,T;\dot{B}_{p,q}^{s})$ and $L^{\rho,r}(0,\infty;\dot{B}_{p,q}^{s})$ will be simplified to $L_{T}^{\rho,r}(\dot{B}_{p,q}^{s})$ and $L^{\rho,r}(\dot{B}_{p,q}^{s})$ respectively. We will also need the following notations: $s_{0}=n/p_{0}+2\alpha/\rho-(2\alpha-1)/m$, $s=n/p+2\alpha/\rho-(2\alpha-1)/m-2\alpha$, $\tilde{s}=s+2\alpha m/\rho$ and $\tilde{\rho}=\rho/(m+1)$, where $p\geq2$, $1<p_{0}\leq p$, $\alpha>1/2$ and $\rho>m+1$ specified later. In addition, subclass of $\dot{B}_{p,q}^{s}$ containing all the solenoidal vectors are denoted by $\dot{B}_{p,q,\sigma}^{s}$.

The main theorem of the paper is posed below.
\begin{thm}\label{thm:1.1}
Assume $m\geq1$, $0<T\leq\infty$, and one of the hypotheses $H_{i}$, $i=0,1,2$ holds. Then there is a small number $0<\eta<1$ such that for all $a\in\dot{B}_{p_{0},r,\sigma}^{s_{0}}$ and $f\in L_{T}^{\tilde{\rho},r}(\dot{B}_{p,\infty}^{\tilde{s}})$ verifying
\begin{equation}\label{cd:af}
\|a\|_{\dot{B}_{p,r}^{s_{0}}}+\|f\|_{L_{T}^{\tilde{\rho},r}(\dot{B}_{p,\infty}^{\tilde{s}})}\leq \eta,
\end{equation}
the generalised Navier-Stokes system \eqref{eqn:gns} has a unique strong solution $(u,\nabla\pi)$ in the class $\big(L_{T}^{\rho,r}(\dot{B}_{p,1,\sigma}^{s+2\alpha})\cap L_{T}^{\tilde{\rho},r}(\dot{B}_{p,\infty}^{\tilde{s}+2\alpha})\big) \times L_{T}^{\tilde{\rho},r}(\dot{B}_{p,\infty}^{\tilde{s}})$ with $\partial_{t}u\in L_{T}^{\tilde{\rho},r}(\dot{B}_{p,\infty}^{\tilde{s}})$. This solution obeys the
following estimates
\begin{align}\label{est:sl}
\nonumber&\quad\|u\|_{L_{T}^{\rho,r}(\dot{B}_{p,1}^{s+2\alpha})}+\|u\|_{L_{T}^{\tilde{\rho},r}(\dot{B}_{p,\infty}^{\tilde{s}+2\alpha})}+\|\partial_{t}u\|_{L_{T}^{\tilde{\rho},r}(\dot{B}_{p,\infty}^{\tilde{s}})}+\|\nabla\pi\|_{L_{T}^{\tilde{\rho},r}(\dot{B}_{p,\infty}^{\tilde{s}})}\\
&\leq C(\|a\|_{\dot{B}_{p_{0},r}^{s_{0}}}+\|f\|_{L_{T}^{\tilde{\rho},r}(\dot{B}_{p,\infty}^{\tilde{s}})})
\end{align}
 for some constant $C>0$ independent of $a$ and $f$.
\end{thm}
\begin{rem}
With the assumption \eqref{cd:af} removed, local well-posedness of (GNS) is still reserved, that is for all $a\in\dot{B}_{p_{0},r,\sigma}^{s_{0}}$ and $f\in L_{T}^{\tilde{\rho},r}(\dot{B}_{p,\infty}^{\tilde{s}})$, there is correspondingly an interval $[0,T')$ with $0<T'<T$ on which the strong solution  $(u,\nabla\pi)$  to \eqref{eqn:gns} is existing and unique in the class  $\big(L_{T'}^{\rho-\delta}(\dot{B}_{p,1,\sigma}^{s+2\alpha})\cap L_{T'}^{\tilde{\rho}-\delta}(\dot{B}_{p,\infty}^{\tilde{s}+2\alpha})\big) \times L_{T'}^{\tilde{\rho}-\delta}(\dot{B}_{p,\infty}^{\tilde{s}})$ with $\partial_{t}u\in L_{T'}^{\tilde{\rho}-\delta,r}(\dot{B}_{p,\infty}^{\tilde{s}})$, where $0<\delta<\tilde{\rho}$ is arbitrary.
\end{rem}

\section{Preliminaries and tool box}

This section is devoted to the investigations on the locally Lipschitz continuity of the high-order power operator $J_{m}$ and the maximal Lorenz regularity of the fractional laplacian $(-\Delta)^{\alpha}$ on homogeneous Besov spaces. Let $C_{0}^{\infty}$ and $\mathcal{S}$ be the families of all smooth functions on $\mathbb{R}^{n}$ having compact supports and rapidly decreasing at the infinity respectively. Let $\mathcal{S}'$ be the dual space of $\mathcal{S}$, whose members are called temperate distributions. For each $1\leq p\leq\infty$, let $L^{p}$ be the common Lebesgue space on $\mathbb{R}^{n}$, either of scalar or vector type, whose norm is denoted by $\|\cdot\|_{p}$.

Given a function $\varphi\in C_{0}^{\infty}$ verifying $0\leq\varphi\leq1$, $\mathrm{supp}\varphi\in \{3/4\leq|\xi|\leq8/3\}$ and
$\sum_{q\in\mathbb{Z}}\varphi(2^{-q}\xi)=1$ for all $\xi\in\mathbb{R}^{n}\setminus\{0\}$. For each $q\in\mathbb{Z}$ and $f\in \mathcal{S}'$, define
$\dot{\Delta}_{q}f=\mathcal{F}^{-1}(\varphi(2^{-q}\xi)\mathcal{F}f)$, where $\mathcal{F}$ and $\mathcal{F}^{-1}$ represent respectively the Fourier transformation and its inverse on $\mathcal{S}'$.
 Then for any $s\in\mathbb{R}$ and $p,r\in[1,\infty]$, we can define the homogeneous Besov space $\dot{B}^{s}_{p,r}$ as follows
\begin{equation*}
  \dot{B}^{s}_{p,r}=\big\{f\in\mathcal{S}_{h}'(\mathbb{R}^{n}):\|f\|_{\dot{B}^{s}_{p,r}}=\|(2^{qs}\|\dot{\Delta}_{q}f\|_{L^{p}})\|_{l^{r}}<\infty\big\},
\end{equation*}
where $\mathcal{S}_{h}'$ is the subset of $\mathcal{S}'$ whose members satisfy $\|\mathcal{F}^{-1}(\psi(\lambda\xi)\mathcal{F}f)\|_{L^{\infty}}\rightarrow0$ as $\lambda\rightarrow\infty$ for all $\psi\in C_{0}^{\infty}$, and $\|(a_{q})\|_{l^{r}}=(\sum_{q\in\mathbb{Z}}|a_{q}|^{r})^{1/r}$ if $1\leq r<\infty$ and $\|(a_{q})\|_{l^{r}}=\sup_{q\in\mathbb{Z}}|a_{q}|$ if $r=\infty$..

Recall that if $0<s<k$ for some $k\in\mathbb{Z}^{+}$, then for all $f\in S_{h}'$, the norm $\|f\|_{\dot{B}_{p,r}^{s}}$ can de characterised by
\begin{align*}
\big(\int_{\mathbb{R}^{n}}\frac{\|\Delta_{y}^{k}f\|_{p}^{r}}{|y|^{sr+n}}dy\big)^{1/r}\;\;\mathrm{if}\;r<\infty\;\;
\mathrm{or}\;\;\sup_{y\in\mathbb{R}^{n}}\frac{\|\Delta_{y}^{k}f\|_{p}}{|y|^{s}}\;\;\mathrm{if}\;r=\infty,
\end{align*}
where $\Delta_{y}f(x)=f(x+y)-f(x)$ and
\begin{align*}
\Delta_{y}^{k}f(x)=\Delta_{y}\Delta_{y}^{k-1}f(x)=\sum_{j=0}^{k}\frac{(-1)^{k-j}k!}{j!(k-j)!}f(x+jy)
\end{align*}
is the $k-$order difference of $f$ (cf. \cite[\S 2.3]{bcd2011}).

The following two lemma can be verified by direct calculations.

\begin{lem}\label{lem:ab}
For all $a,b\in\mathbb{R}^{n}$,
\begin{align*}
|J_{m}(a)-J_{m}(b)|
\leq\left\{\begin{array}{l}
                                  m(|a|^{m-1}+|b|^{m-1})|a-b|, \quad m>1, \\
                                  6|a-b|^{m},\quad 0<m\leq1.
                                \end{array}\right.
\end{align*}
\end{lem}
\noindent {\it Proof} : Since the case $m=1$ is evident, we divide the proof into two cases: $0<m<1$ and $1<m<\infty$. In the case $0<m<1$, without loss of generality, we may assume $|a|\geq|b|>0$. If $|a|\geq2|b|$, then $|a-b|\geq|a|/2$. So we can let $e_{a}=a/|a|$ and $c=b/|a|$ to deduce that
 \begin{align*}
&\quad|J_{m}(a)-J_{m}(b)|=|a|^{m}|e_{a}-J_{m}(c)|\\
&\leq(1+|c|^{m})|a|^{m}\leq (1+2^{-m})2^{m}|a-b|^{m}\leq3|a-b|^{m}.
\end{align*}
If $|b|\leq|a|<2|b|$, then we have $|a-b|\leq3|b|$. Consequently,
 \begin{align*}
&\quad|J_{m}(a)-J_{m}(b)|\leq|b|^{m-1}|a-b|+(|b|^{m-1}-|a|^{m-1})|a|\\
&\leq|b|^{m-1}|a-b|+(|b|^{m-1}|a|-|b|^{m}))\leq2|b|^{m-1}|a-b|\\
&\leq 2\cdot3^{1-m}|a-b|^{m-1}|a-b|\leq6|a-b|^{m}.
\end{align*}

In the case $1<m<\infty$, the vector-valued function $J_{m}(a)$ is of $C^{1}$ type on $\mathbb{R}^{n}$ with $\nabla J_{m}(0)=0$ and $\nabla J_{m}(a)=|a|^{m-1}E+(m-1)|a|^{m-3}a\otimes a$ for $a\neq0$, where $E$ is the $n-$order identity, and $a\otimes a$ is the tensor product of $a$ and $a$. Hence,
\begin{align}\label{ine:mab}
\nonumber&\quad|J_{m}(a)-J_{m}(b)|=\big|\int_{0}^{1}(a-b)\cdot\nabla J_{m}((1-\lambda)a+\lambda b)d\lambda\big|\\
&\leq m\big|\int_{0}^{1}|(1-\lambda)a+\lambda b|^{m-1}d\lambda\cdot|a-b|\leq m(|a|^{m-1}+|b|^{m-1})|a-b|.
\end{align}
Thus the lemma has been proved. \hfill $\Box$

\begin{prop}\label{prop:pm1}
Assume that $0<s<m\leq1$ and $1/m\leq p,r\leq\infty$, then for all $f\in \dot{B}^{s/m}_{mp,mr}$, we have $J_{m}(f)\in \dot{B}^{s}_{p,r}$, and
\begin{align}\label{est:pm1}
\|J_{m}(f)\|_{\dot{B}^{s}_{p,r}}\leq C\|f\|_{\dot{B}^{s/m}_{mp,mr}}^{m}
\end{align}
for some constant $C>0$ independent of $f$.
\end{prop}
\noindent {\it Proof} : The case $m=1$ is trivial, so we only consider the case $m<1$. We also assume that $r<\infty$, the case $r=\infty$ can be dealt with in a similar way.

Since $\dot{B}^{s/m}_{mp,mr}\subseteq L_{\mathrm{loc}}^{mp}$, we first have $J_{m}(f)\in L_{\mathrm{loc}}^{p}$. Then by pursuing the same line of reasoning as in the proof of \cite[Lemma 2.12]{bcd2011}, we can check that $J_{m}(f)\in S_{h}'$.

For all $x,y\in\mathbb{R}^{n}$, by invoking Lemma \ref{lem:ab}, we obtain
\begin{align*}
|\Delta_{y}J_{m}(f)(x)|&=|J_{m}(f(x+y))-J_{m}(f(x))|
\leq6|\Delta_{y}f(x)|^{m}.
\end{align*}
Consequently,
\begin{align*}
\big(\int_{\mathbb{R}^{n}}\frac{\|\Delta_{y}J_{m}(f)\|_{p}^{r}}{|y|^{sr}}\frac{dy}{|y|^{n}}\big)^{1/r}
&\leq6\big(\int_{\mathbb{R}^{n}}\frac{\|\Delta_{y}f\|_{mp}^{mr}}{|y|^{sr}}\frac{dy}{|y|^{n}}\big)^{1/r},
\end{align*}
which, together with the equivalent characterization of the norm of $\dot{B}^{s}_{p,r}$, yields (\ref{est:pm1}) immediately.  \hfill $\Box$

 \begin{lem}[\cite{bcd2011}, \S 2.6 or \cite{mwz2012}, \S 3.1 ]\label{est:pro}
Let $s>0$, $p,p_{i},q_{i},r\in[1,\infty]$, $\delta_{i}>0$, $i=1,2$ such that $1/p=1/p_{1}+1/p_{2}=1/q_{1}+1/q_{2}$ and $s<n/p$ or $s=n/p$ and $r=1$,
then we have
\begin{align}\label{est:pro1}
 \|fg\|_{\dot{B}^{s}_{p,r}}\leq C(\|f\|_{\dot{B}^{s+\delta_{1}}_{p_{1},r}}\|g\|_{\dot{B}^{-\delta_{1}}_{p_{2},\infty}}+\|f\|_{\dot{B}^{-\delta_{2}}_{q_{1},\infty}}\|g\|_{\dot{B}^{s+\delta_{2}}_{q_{2},r}})
 \end{align}
or
 \begin{align}\label{est:pro2}
\|fg\|_{\dot{B}^{s}_{p,r}}\leq C(\|f\|_{\dot{B}^{s}_{p_{1},r}}\|g\|_{p_{2}}+\|f\|_{q_{1}}\|g\|_{\dot{B}^{s}_{q_{2},r}}).
\end{align}
\end{lem}

\begin{prop}
Suppose that $m\geq1$, $1<p<\infty$, $1\leq r\leq\min\{2,p\}$ and $0<s<\min\{m,n/p\}$, then it holds
\begin{equation}\label{est:pm11}
\|J_{m}(f)\|_{\dot{B}_{p,r}^{s}}\leq C\|f\|_{\dot{B}_{p,r}^{s/m+(1-1/m)n/p}}^{m},
\end{equation}
where $C=C(s,p,n,m)>0$.
\end{prop}
{\it Proof} : We only deal with $m>1$, the case $m=1$ is trivial.  By pursuing the same arguments as in \cite[\S 1.5]{mwz2012} with slight modification, we have
\begin{equation}\label{est:pm2}
\|J_{m}(f)\|_{\dot{B}_{p,r}^{s}}\leq C\|f\|_{p_{1}}^{m-1}\|f\|_{\dot{B}_{p_{2},r}^{s}},
\end{equation}
 where $1\leq r\leq\infty$, $\max\{1,p(m-1)\}< p_{1}<\infty$, $p<p_{2}<\infty$ verifying $1/p=(m-1)/p_{1}+1/p_{2}$.

Take $p_{1},p_{2}$ such that
\begin{align*}
\frac{n}{p_{1}}=\frac{1}{m}\big(\frac{n}{p}-s\big)\;\;\textrm{and}\;\;\frac{n}{p_{2}}=\big(1-\frac{1}{m}\big)s+\frac{n}{mp},
\end{align*}
then under the assumptions on $s,p,r,m$, we have $\max\{1,m-1\}p< p_{1}<\infty$, $p<p_{2}<\infty$ and $s/m+(1-1/m)n/p=n/p-n/p_{1}=s+n/p-n/p_{2}$. Consequently all the following embedding relations come true.
\begin{align*}
\dot{B}^{s/m+(1-1/m)n/p}_{p,r}\hookrightarrow\dot{B}^{0}_{p_{1},r}\cap \dot{B}^{s}_{p_{2},r},
\end{align*}
and
\begin{align*}
\dot{B}^{0}_{p_{1},r}\hookrightarrow\left\{\begin{array}{ll}
                              \dot{B}^{0}_{p_{1},p}\hookrightarrow\dot{B}^{0}_{p_{1},p_{1}}, & \textrm{if}\;p<p_{1}\leq2, \\
                              \dot{B}^{0}_{p_{1},p}\hookrightarrow\dot{B}^{0}_{p_{1},2}, & \textrm{if}\;p<2<p_{1}, \\
                               \dot{B}^{0}_{p_{1},2}, & \textrm{if}\;2\leq p<p_{1}\\
                             \end{array}
                             \right\}\hookrightarrow L^{p_{1}}.
\end{align*}
Please see \cite[\S 2.5]{bcd2011} for references. Putting these embedding relations into (\ref{est:pm2}), we obtain the desired estimate (\ref{est:pm11}).  \hfill $\Box$

\begin{prop}\label{thm:pm}
Assume that $1<p<\infty$,
\begin{align}\label{cd:sm}
\left\{\begin{array}{ll}
0<s<\min\{m-1,(m-1)^{2}n/p\},\;1/(m-1)\leq r\leq\infty,\\
\qquad\qquad\qquad\qquad\qquad\qquad\qquad\qquad\qquad\textrm{if}\;\; 1<m<2,\\
           0<s<\min\{m-1, n/p\},\;1\leq r\leq\min\{2,p\},\;\;\textrm{if}\;\;m\geq 2.
           \end{array}
\right.
\end{align}
Then there exists a constant $C>0$ such that
\begin{align}\label{ine:mf21}
\nonumber&\quad\|J_{m}(f)-J_{m}(g)\|_{\dot{B}_{p,r}^{s}}\\
&\leq C(\|f\|_{\dot{B}^{s/m+(1-1/m)n/p}_{p,r_{0}}}^{m-1}+\|g\|_{\dot{B}^{s/m+(1-1/m)n/p}_{p,r_{0}}}^{m-1})\|f-g\|_{\dot{B}^{s/m+(1-1/m)n/p}_{p,r_{0}}},
\end{align}
where $r_{0}=1$ if $1<m<2$, and $r_{0}=r$ if $m\geq 2$.
\end{prop}
\noindent {\it Proof} : For any $x,y\in\mathbb{R}^{n}$, by virtue of \eqref{ine:mab}, we can derive that
\begin{align*}
J_{m}(f(x))-J_{m}(g(x))=\int_{0}^{1}J'_{m}((1-\lambda)f(x)+\lambda g(x))d\lambda\cdot(f(x)-g(x)).
\end{align*}
Then by employing (\ref{est:pro2}), we obtain
\begin{align}\label{ine:mfg}
\nonumber\|J_{m}(f)-J_{m}(g)\|_{\dot{B}_{p,r}^{s}}
&\leq C\big(\int_{0}^{1}\|J'_{m}((1-\lambda)f+\lambda g)\|_{\dot{B}^{s}_{p_{1},r}}d\lambda\cdot\|f-g\|_{p_{2}}\\
&\quad+\int_{0}^{1}\|J'_{m}((1-\lambda)f+\lambda g)\|_{q_{1}}d\lambda\cdot\|f-g\|_{\dot{B}^{s}_{q_{2},r}}\big),
\end{align}
where $1/p=1/p_{1}+1/p_{2}=1/q_{1}+1/q_{2}$ and $r\in[1,\infty]$ is free.

In the case $1<m<2$, we can use (\ref{est:pm1}) to deduce that
\begin{align}\label{ine:mf2}
\nonumber\|J_{m}(f)-J_{m}(g)\|_{\dot{B}_{p,r}^{s}}
\nonumber&\leq C\int_{0}^{1}\|(1-\lambda)f+\lambda g\|_{\dot{B}^{s/(m-1)}_{(m-1)p_{1},(m-1)r}}^{m-1}d\lambda\cdot\|f-g\|_{p_{2}}\\
&\quad+C\int_{0}^{1}\|(1-\lambda)f+\lambda g\|_{(m-1)q_{1}}^{m-1}d\lambda\cdot\|f-g\|_{\dot{B}^{s}_{q_{2},r}}\\
\nonumber&\leq C(\|f\|_{\dot{B}^{s/(m-1)}_{(m-1)p_{1},1}}^{m-1}+\|g\|_{\dot{B}^{s/(m-1)}_{(m-1)p_{1},1}}^{m-1})\|f-g\|_{p_{2}}\\
\nonumber&\quad+C(\|f\|_{(m-1)q_{1}}^{m-1}+\|g\|_{(m-1)q_{1}}^{m-1})\|f-g\|_{\dot{B}^{s}_{q_{2},1}}.
\end{align}
Take $p/(m-1)\leq p_{1},q_{1}<\infty$  such that
\begin{align}
\label{cd:pi}&\frac{n}{p_{1}}=\big(1-\frac{1}{m}\big)\frac{n}{p}+\frac{s}{m},\;\;\frac{n}{p_{2}}=\frac{1}{m}\big(\frac{n}{p}-s\big),\;\;\textrm{and}\\
\label{cd:qi}&\frac{n}{q_{1}}=\big(1-\frac{1}{m}\big)\big(\frac{n}{p}-s\big),\;\;\frac{n}{q_{2}}=\frac{n}{mp}+\big(1-\frac{1}{m}\big)s.
\end{align}
Under this setting, we have
$p\leq(m-1)p_{1}\leq p_{2}$, $p\leq q_{2}\leq (m-1)q_{1}$,
and then
\begin{align*}
\dot{B}^{s/m+(1-1/m)n/p}_{p,1}\hookrightarrow\dot{B}^{s/(m-1)}_{(m-1)p_{1},1}\hookrightarrow L^{p_{2}},\;\;\dot{B}^{s/m+(1-1/m)n/p}_{p,1}\hookrightarrow\dot{B}^{s}_{q_{2},1}\hookrightarrow L^{(m-1)q_{1}}.
\end{align*}
Inserting these embedding relations into (\ref{ine:mf2}), we obtain (\ref{ine:mf21}).

In the case $m\geq 2$, we can also take $p\leq p_{i},q_{i}<\infty$, $i=1,2$ as in (\ref{cd:pi}) and (\ref{cd:qi}) respectively.
Under this setting,
it is easy to see that $s<n/p_{1}$ if and only if $s<n/p$. Consequently,
by virtue of (\ref{est:pm2}) with $m$ replaced by $m-1$ and the following embedding
\begin{align*}
 \dot{B}^{s/m+(1-1/m)n/p}_{p,r}\hookrightarrow\dot{B}^{s/(m-1)+(1-1/(m-1))n/p_{1}}_{p_{1},r},
\end{align*}
we can derive from (\ref{ine:mfg}) that
\begin{align*}
&\quad\|J_{m}(f)-J_{m}(g)\|_{\dot{B}_{p,r}^{s}}\\
&\leq C\int_{0}^{1}\big\|(1-\lambda)f+\lambda g\big\|_{\dot{B}^{s/(m-1)+(1-1/(m-1))n/p_{1}}_{p_{1},r}}^{m-1}d\lambda\cdot\|f-g\|_{p_{2}}\\
&\quad+C\int_{0}^{1}\|(1-\lambda)f+\lambda g\|_{(m-1)q_{1}}^{m-1}d\lambda\cdot\|f-g\|_{\dot{B}^{s}_{q_{2},r}}\\
&\leq C(\|f\|_{\dot{B}^{s/m+(1-1/m)n/p}_{p,1}}^{m-1}+\|g\|_{\dot{B}^{s/m+(1-1/m)n/p}_{p,r}}^{m-1})\|f-g\|_{\dot{B}^{s/m+(1-1/m)n/p}_{p,r}}.
\end{align*}
Thus the proof has been completed. \hfill $\Box$

\begin{rem}
Derivation of \eqref{ine:mf21} for $1<m<2$ relies on the estimate \eqref{est:pm1}. This is the reason why condition \eqref{cd:sm} for $1<m<2$ is much stricter than that for $m\geq2$. We wonder whether or not the estimate \eqref{ine:mf21} keeps true under the weaker assumption $0<s<\min\{m-1,n/p\}$.
\end{rem}

Given $\alpha>0$, the fractional power of the negative Laplacian $(-\Delta)^{\pm\alpha}$ is defined by $\mathcal{F}^{-1}|\xi|^{\pm2\alpha}\mathcal{F}$. Since the function $|\xi|^{\pm2\alpha}$ is $\pm2\alpha-$order homogeneous on $\mathbb{R}^{n}\setminus\{0\}$, $(-\Delta)^{\alpha}$ is an isomorphism from $\dot{B}^{s+2\alpha}_{p,r}$ onto $\dot{B}^{s}_{p,r}$ provided $s+2\alpha<n/p$, or $s+2\alpha=n/p$ and $r=1$. Define the generalised heat semigroup $e^{-t (-\Delta)^{\alpha}}=\mathcal{F}^{-1}e^{-t|\xi|^{2\alpha}}\mathcal{F}$ for $t\geq0$.
 Next proposition gives an estimate for $e^{-t (-\Delta)^{\alpha}}$ on Lorentz-Besov spaces (cf. \cite{ks2019,zcz2024}).
\begin{prop}\label{prop:1.1}
Let $0\leq \gamma\leq\alpha$, $1\leq r\leq\infty$, and let $s_{0},s\in \mathbb{R}$, $1<\rho<\infty$ and $1\leq p_{0}\leq p\leq\infty$ such that
 \begin{align}\label{cd:sp1}
s_{0}-\frac{n}{p_{0}}-2\alpha<s-\frac{n}{p}<s_{0}-\frac{n}{p_{0}}\;\;\mathrm{and}\;\;s-\frac{n}{p}-\frac{2\alpha}{\rho}=s_{0}-\frac{n}{p_{0}}-2\alpha.
 \end{align}
 If $a\in \dot{B}_{p_{0},r}^{s_{0}}$, then $(-\Delta)^{\gamma}e^{-t (-\Delta)^{\alpha}}a\in L^{\rho,r}(\dot{B}_{p,1}^{s+2(\alpha-\gamma)})$, and
\begin{equation}\label{ine:exp1}
 \|(-\Delta)^{\gamma}e^{-t (-\Delta)^{\alpha}}a\|_{L^{\rho,r}(\dot{B}_{p,1}^{s+2(\alpha-\gamma)})}\leq C\|a\|_{\dot{B}_{p_{0},r}^{s_{0}}},
\end{equation}
where the constant $C>0$ is unconnected with $a$ and $\gamma$.
\end{prop}

The following proposition reveals the maximal Lorentz regularity of $(-\Delta)^{\alpha}$. This is an generalisation of that in \cite{ks2019}, and was proved in \cite{zcz2024}.
\begin{prop}\label{prop:1.1a}
Under hypothesis (\ref{cd:sp1}) with $1<p_{0}\leq p<\infty$ and $1\leq q,r\leq\infty$, for every $0<T\leq\infty$, any $a\in \dot{B}_{p_{0},r}^{s_{0}}$ and $f\in L_{T}^{\rho,r}(\dot{B}_{p,q}^{s})$, Cauchy problem
\begin{equation}\label{eq:evo}
 \frac{du}{dt}+(-\Delta)^{\alpha}u=f(t),\;t>0,\; u(0)=a
\end{equation}
has a unique solution $u$ such that both $u'$ and $(-\Delta)^{\alpha}u$ lie in the class $L_{T}^{\rho,r}(\dot{B}_{p,q}^{s})$,  and
\begin{equation*}
  \|u'\|_{L_{T}^{\rho,r}(\dot{B}_{p,q}^{s})}+  \|(-\Delta)^{\alpha}u\|_{L_{T}^{\rho,r}(\dot{B}_{p,q}^{s})}\leq C\big(\|a\|_{\dot{B}_{p_{0},r}^{s_{0}}}+ \|f\|_{L_{T}^{\rho,r}(\dot{B}_{p,q}^{s})}\big)
\end{equation*}
for some constant $C>0$ free from  $a$, $T$ and $f$.
\end{prop}

Evidently, each solution of (\ref{eq:evo}) has the following representation
\begin{align}\label{rep:u}
u(t)=e^{-t(-\Delta)^{\alpha}}a+\int_{0}^{t}e^{-(t-\tau)(-\Delta)^{\alpha}}f(\tau)d\tau=:a_{L}(t)+(Sf)(t),
\end{align}
where the function $Sf$ solves problem (\ref{eq:evo}) with $u(0)=0$, and assumes more regularity than $f$ does, that is (see \cite{ks2019,zcz2024} for references).
\begin{prop}\label{prop:1.2}
Assume that $s,\tilde{s}\in\mathbb{R}$, $1\leq r\leq\infty$, $1<\tilde{p}\leq p<\infty$ and $1<\tilde{\rho}<\rho<\infty$ verifying
\begin{align}\label{cd:sp2}
\tilde{s}-\frac{n}{\tilde{p}}-2\alpha<s-\frac{n}{p}<\tilde{s}-\frac{n}{\tilde{p}}\;\;\mathrm{and}\;\;s-\frac{n}{p}-\frac{2\alpha}{\rho}=\tilde{s}-\frac{n}{\tilde{p}}-\frac{2\alpha}{\tilde{\rho}},
 \end{align}
 then for each $0\leq\gamma\leq\alpha$, it holds
 \begin{equation}\label{est:f1}
  \|(-\Delta)^{\gamma}(Sf)\|_{L_{T}^{\rho,r}(\dot{B}_{p,1}^{s+2(\alpha-\gamma)})}\leq C\|f\|_{L^{\tilde{\rho},r}(0,T;\dot{B}_{\tilde{p},\infty}^{\tilde{s}})},
 \end{equation}
 where the constant $C>0$ is independent of $\gamma$, $f$ and $0<T\leq\infty$.
\end{prop}

The following lemma can be verified by means of interpolation (refer to \cite{ks2019} for $k=2$).
\begin{lem}\label{lem:pfg}
Let $1<\rho\leq\rho_{i}<\infty$, $i=1,2,\cdots,k$ such that $\sum_{i=1}^{k}1/\rho_{i}=1/\rho$ and $1\leq r\leq \infty$, then we have
\begin{equation*}
\big\|\prod_{i=1}^{k}f_{i}\big\|_{L^{\rho,r}(\mathbb{R})}\leq C\prod_{i=1}^{k}\|f_{i}\|_{L^{\rho_{i},r}(\mathbb{R})},
\end{equation*}
where the constant $C>0$ is not related to $f_{i}$.
\end{lem}
\begin{prop}[\cite{zcz2024}]\label{lem:pm}
Let $1<\rho<\infty$, $1\leq r\leq\infty$ and $m\geq1$, then there exists a constant $C>0$ such that
\begin{equation*}
\frac{ 1}{C}\|\varphi\|_{L^{m\rho,mr}(\mathbb{R})}^{m}\leq\||\varphi|^{m}\|_{L^{\rho,r}(\mathbb{R})}\leq C\|\varphi\|_{L^{m\rho,mr}(\mathbb{R})}^{m}.
\end{equation*}
\end{prop}

Define $\Phi_{m}(u,v)=-S(P(J_{m}(u)\cdot\nabla v))$.
Let $m+1<\rho<\infty$, and take
\begin{align}\label{val:s}
 s=\frac{n}{p}+\frac{2\alpha}{\rho}-\frac{2\alpha-1}{m}-2\alpha,\;\tilde{s}=s+\frac{2m\alpha}{\rho}\;\textrm{and}\;\tilde{\rho}=\frac{\rho}{m+1}.
\end{align}
Then thanks to Proposition \ref{prop:1.2}, we have
\begin{align}\label{est:phu}
&\quad\|\Phi_{m}(u,v)\|_{L_{T}^{\rho,r}(\dot{B}_{p,1}^{s+2\alpha})}+\|(\Phi_{m}(u,v))'\|_{L_{T}^{\rho,r}(\dot{B}_{p,1}^{s})}\\
\nonumber&\leq C\|J_{m}(u)\cdot\nabla v\|_{L_{T}^{\tilde{\rho},r}(\dot{B}_{p,\infty}^{\tilde{s}})},
\end{align}
where $C>0$ is independent of $u,v$ and $T$.

\section{Proofs of the main theorems}

We first give the definition for the strong solution to the (GBS). Given $s,s_{0}\in \mathbb{R}$, $p\in(1,\infty)$, $r,q\in[1,\infty)$ and $a\in \dot{B}_{p,r,\sigma}^{s_{0}}$. A function couple $(u,\nabla\pi)$ is said to be the strong solution of (\ref{eqn:gns}) on the interval $[0,T]$, if
 \begin{align*}
u\in L^{1}(0,T;\dot{B}_{p,q,\sigma}^{s+2\alpha})\cap C([0,T);\dot{B}_{p,r,\sigma}^{s_{0}}+\dot{B}_{p,q,\sigma}^{s}),\; u',\nabla\pi\in L^{1}(0,T;\dot{B}_{p,q}^{s})
 \end{align*}
with $u(0)=a$, and  the equation of (\ref{eqn:gns}) is verified by $u(t)$ and $\nabla\pi(t)$ in $\dot{B}_{p,q}^{s}$ for a.e. $t\in[0,T]$.

By performing the Helmholtz projection $\mathcal{P}$ on both sides of \eqref{eqn:gns}, an abstract evolution equation,
\begin{equation}\label{eqn:gns1}
\partial_{t}u+\mu(-\Delta)^{\alpha}u=-\mathcal{P}(J_{m}(u)\cdot\nabla u)+\mathcal{P}f(t),\;\;t>0, \;\; u(0)=a.
\end{equation}
is derived, where commutativity of $\mathcal{P}$ and $(-\Delta)^{\alpha}$ is taken into account. Evidently, for any strong solution $(u,\nabla\pi)$ of (\ref{eqn:gns}), the first function $u$ solves \eqref{eqn:gns1} in the space $\dot{B}_{p,q}^{s}$  a.e. on $[0,T]$. Hence, to investigate the strong solvability of (\ref{eqn:gns}), it is naturally to search for the strong solution of \eqref{eqn:gns1}, where the nonlinear term $\mathcal{P}(J_{m}(u)\cdot\nabla v)$ will be treated firstly. Note that the operator $\mathcal{P}$ is bounded on Besov spaces, so in the following arguments we will omit it, and only take $J_{m}(u)\cdot\nabla v$ into account.

\begin{prop}\label{prop:3.1}
Assume that

\noindent$H_{0}$: $m=1$, $1<p<\infty$, $2<\rho<\infty$ and
\begin{align*}
s=1+\frac{n}{p}+\frac{2\alpha}{\rho}-4\alpha,\;\;\tilde{s}=s+\frac{2\alpha}{\rho}
\end{align*}
such that
  \begin{align}\label{hy:a0}
\frac{1}{2}<\alpha<1+\frac{n}{2p}\;\;\textrm{and}\;\;
2\alpha-1-\frac{n}{2p}<\frac{2\alpha}{\rho}<2\alpha-1.
\end{align}
Then for all $u,v\in L_{T}^{\rho,r}(\dot{B}_{p,1}^{s+2\alpha})$, the following inequality
\begin{align}\label{est:m1}
\|u\cdot\nabla v\|_{L^{\rho/2,r}(0,T;\dot{B}_{p,\infty}^{\tilde{s}})}\leq C\|u\|_{L_{T}^{\rho,r}(\dot{B}_{p,1}^{s+2\alpha})}\|v\|_{L_{T}^{\rho,r}(\dot{B}_{p,1}^{s+2\alpha})}
\end{align}
holds true.
\end{prop}
\noindent{\it Proof} : For a.e. $t\in(0,T)$, it can be deduced from \eqref{est:pro1} that
\begin{align*}
\|u(t)\cdot\nabla v(t)\|_{\dot{B}_{p,\infty}^{\tilde{s}}}&\leq\|u(t)\otimes v(t)\|_{\dot{B}_{p,\infty}^{1+\tilde{s}}}\\
&\leq C\big(\|u(t)\|_{\dot{B}_{\infty,\infty}^{-\delta_{0}}}\|v(t)\|_{\dot{B}_{p,\infty}^{1+\tilde{s}+\delta_{0}}}+\|v(t)\|_{\dot{B}_{\infty,\infty}^{-\delta_{0}}}\|u(t)\|_{\dot{B}_{p,\infty}^{1+\tilde{s}+\delta_{0}}}\big)\\
&\leq C\big(\|u(t)\|_{\dot{B}_{p,\infty}^{n/p-\delta_{0}}}\|v(t)\|_{\dot{B}_{p,\infty}^{1+\tilde{s}+\delta_{0}}}+\|v(t)\|_{\dot{B}_{p,\infty}^{n/p-\delta_{0}}}\|u(t)\|_{\dot{B}_{p,\infty}^{1+\tilde{s}+\delta_{0}}}\big)\\
&\leq C\|u(t)\|_{\dot{B}_{p,1}^{s+2\alpha}}\|v(t)\|_{\dot{B}_{p,1}^{s+2\alpha}}.
\end{align*}
Then by invoking Lemma \ref{lem:pfg} with $k=2$ and $f_{1},f_{2}$ replaced by $\|u(t)\|_{\dot{B}_{p,1}^{s+2\alpha}},\|v(t)\|_{\dot{B}_{p,1}^{s+2\alpha}}$ respectively, the desired inequality \eqref{est:m1} is easily reached. \hfill $\Box$

\begin{prop}\label{prop:3.2}
Let $m>1$, $m+1<\rho<\infty$ and $s,\tilde{s},\rho,\tilde{\rho}$ be as in \eqref{val:s}, then
under one of the following assumptions:

\noindent$H_{1}$: $1<m<2$,
\begin{align}\label{hy:p1}
\max\big\{\frac{(m-1)n}{3-m},\frac{4-m}{3-m}\big\}<p<\frac{n}{m-1},
\end{align}
\begin{align}\label{hy:a1}
\frac{1}{2}+\frac{m(2-m)n}{2(3-m)p}<\alpha<\frac{m+1}{2}+\frac{mn}{2p},
\end{align}
and
\begin{align}\label{hy:r1}
\frac{2\alpha-1}{m}-\frac{n}{(m+1)p}<\frac{2\alpha}{\rho}<\frac{2\alpha-1}{m}-\frac{(2-m)n}{(3-m)p},
\end{align}
\noindent$H_{2}$: $2\leq m<\infty$,  $n\leq p<2n$,
\begin{align}\label{hy:a2}
\frac{1}{2}<\alpha<\frac{1}{2}+\frac{mn}{p},
\end{align}
and
\begin{align}\label{hy:r2}
\frac{2\alpha-1}{m}+\frac{1}{m+1}\big(1-\frac{2n}{p}\big)
<\frac{2\alpha}{\rho}<\frac{2\alpha-1}{m},
\end{align}
it holds that
\begin{equation}\label{est:uv1}
\|J_{m}(u)\cdot\nabla v\|_{L^{\tilde{\rho},r}(0,T;\dot{B}_{p,\infty}^{\tilde{s}})}\leq C\|u\|_{L_{T}^{\rho,r}(\dot{B}_{p,\infty}^{s+2\alpha})}^{m}\|v\|_{L_{T}^{\rho,r}(\dot{B}_{p,\infty}^{s+2\alpha})},
\end{equation}
and
\begin{align}\label{est:uv2}
\nonumber&\quad\|(J_{m}(u_{1})-J_{m}(u_{2}))\cdot\nabla v\|_{L_{T}^{\tilde{\rho},r}(\dot{B}_{p,\infty}^{\tilde{s}})}\\
&\leq C(\|u_{1}\|_{L_{T}^{\rho,r}(\dot{B}_{p,1}^{s+2\alpha})}^{m-1}+\|u_{2}\|_{L_{T}^{\rho,r}(\dot{B}_{p,1}^{s+2\alpha})}^{m-1})
\|u_{1}-u_{2}\|_{L_{T}^{\rho,r}(\dot{B}_{p,1}^{s+2\alpha})}\|v\|_{L_{T}^{\rho,r}(\dot{B}_{p,1}^{s+2\alpha})}.
\end{align}
\end{prop}
\noindent{\it Proof} : We only prove \eqref{est:uv2}, and \eqref{est:uv1} is a consequence of \eqref{est:uv2} where $u_{1}=u$ and $u_{2}=0$. For the sake of convenience, we ignore the dependence of $u_{i},v$ on $t$ temporarily and assume that $v\in\mathcal{S}_{0}^{n}$, i.e. $v\in \mathcal{S}^{n}$, and its Fourier transform vanishes near the original.

Take any $\varphi\in\mathcal{S}_{0}$ and consider the dual product $\langle J_{m}(u_{1})-J_{m}(u_{2}))\cdot\nabla v,\varphi\rangle$. By the generalised divergence formula, we have
\begin{align*}
&\quad\langle J_{m}(u_{1})-J_{m}(u_{2}))\cdot\nabla v,\varphi\rangle\\
&=\langle\textrm{div}((J_{m}(u_{1})-J_{m}(u_{2}))\otimes v),\varphi\rangle-\langle\textrm{div}(J_{m}(u_{1})-J_{m}(u_{2})) v,\varphi\rangle\\
&=\langle\textrm{div}((J_{m}(u_{1})-J_{m}(u_{2}))\otimes v),\varphi\rangle+\langle J_{m}(u_{1})-J_{m}(u_{2}),\nabla(v\cdot\varphi)\rangle\\
&=:I_{1}+I_{2}.
\end{align*}

When $1<m<2$, we take $p_{1}=(3-m)p$, $p_{2}=(3-m)p/(2-m)$,
\begin{align*}
\delta_{1}=2\alpha-1-\frac{2\alpha m}{\rho},\;\;\delta_{2}=\frac{n}{p_{1}}-s-2\alpha\;\;\textrm{and}\;\;
\sigma_{1}=\frac{n}{p_{1}}-\frac{\delta_{1}}{m}.
\end{align*}
Under hypotheses $H_{1}$, it holds that $\delta_{i}>0$, $i=1,2$,
\begin{align*}
-\min\big\{1,\frac{n}{p_{1}}\big\}<\tilde{s}<0,
\end{align*}
and
\begin{align*}
0<\frac{n}{p_{1}}-\delta_{1}=\tilde{s}+1+\delta_{2}<\frac{(m-1)^{2}n}{p_{1}}<m-1.
\end{align*}
Then by the dual principle for $\dot{B}_{p,\infty}^{\tilde{s}}\times\dot{B}_{p',1}^{-\tilde{s}}$, where $p'=p/(p-1)$ is conjugate number of $p$, together with inequalities \eqref{est:pro1} and \eqref{ine:mf21}, we have
\begin{align*}
|I_{1}|&\leq\|\textrm{div}((J_{m}(u_{1})-J_{m}(u_{2}))\otimes v)\|_{\dot{B}_{p,\infty}^{\tilde{s}}}\|\varphi\|_{\dot{B}_{p',1}^{-\tilde{s} }}\\
&\leq C\|(J_{m}(u_{1})-J_{m}(u_{2}))\otimes v\|_{\dot{B}_{p,\infty}^{\tilde{s}+1}}\|\varphi\|_{\dot{B}_{p',1}^{-\tilde{s} }}\\
&\leq C\big(\|J_{m}(u_{1})-J_{m}(u_{2})\|_{\dot{B}_{p_{1},\infty}^{n/p_{1}-\delta_{1}}}\| v\|_{\dot{B}_{p,\infty}^{\tilde{s}+1+\delta_{1}}}\\
&\quad+\|J_{m}(u_{1})-J_{m}(u_{2})\|_{\dot{B}_{p_{1},\infty}^{\tilde{s}+1+\delta_{2}}}\|v\|_{\dot{B}_{p_{2},\infty}^{-\delta_{2}}}\big)\|\varphi\|_{\dot{B}_{p',1}^{-\tilde{s} }}\\
&\leq C\big[\big(\|u_{1}\|_{\dot{B}_{p_{1},1}^{\sigma_{1}}}^{m-1}+\|u_{2}\|_{\dot{B}_{p_{1},1}^{\sigma_{1}}}^{m-1}\big)\|u_{1}-u_{2}\|_{\dot{B}_{p_{1},1}^{\sigma_{1}}}\|v\|_{\dot{B}_{p,1}^{\tilde{s}+\delta_{1}+1}}\\
&\quad+\big(\|u_{1}\|_{\dot{B}_{p_{1},1}^{\sigma_{1}}}^{m-1}+\|u_{2}\|_{\dot{B}_{p_{1},1}^{\sigma_{1}}}^{m-1}\big)\|u_{1}-u_{2}\|_{\dot{B}_{p_{1},1}^{\sigma_{1}}}\|v\|_{\dot{B}_{p_{2},\infty}^{-\delta_{2}}}\big]\|\varphi\|_{\dot{B}_{p',1}^{-\tilde{s}}}\\
&\leq C\big[\big(\|u_{1}\|_{\dot{B}_{p,1}^{\sigma_{1}+n/p_{2}}}^{m-1}+\|u_{2}\|_{\dot{B}_{p,1}^{\sigma_{1}+n/p_{2}}}^{m-1}\big)\|u_{1}-u_{2}\|_{\dot{B}_{p,1}^{\sigma_{1}+n/p_{2}}}\|v\|_{\dot{B}_{p,1}^{\tilde{s}+\delta_{1}+1}}\\
&\quad+\big(\|u_{1}\|_{\dot{B}_{p,1}^{\sigma_{1}+n/p_{2}}}^{m-1}+\|u_{2}\|_{\dot{B}_{p,1}^{\sigma_{1}+n/p_{2}}}^{m-1}\big)\|u_{1}-u_{2}\|_{\dot{B}_{p,1}^{\sigma_{1}+n/p_{2}}}\|v\|_{\dot{B}_{p,\infty}^{n/p_{1}-\delta_{2}}}\big]\|\varphi\|_{\dot{B}_{p',1}^{-\tilde{s}}}\\
&\leq C(\|u_{1}\|_{\dot{B}_{p,1}^{s+2\alpha}}^{m-1}+\|u_{2}\|_{\dot{B}_{p,1}^{s+2\alpha}}^{m-1})\|u_{1}-u_{2}\|_{\dot{B}_{p,1}^{s+2\alpha}}\|v\|_{\dot{B}_{p,1}^{s+2\alpha}}\|\varphi\|_{\dot{B}_{p',1}^{-\tilde{s} }}.
\end{align*}
As for $I_{2}$, we take $p'<q_{1}<\infty$ and $p<q_{2}<\infty$ such that
\begin{align*}
\frac{1}{q_{1}}=\frac{1-\eta_{1}}{p'}\;\;\textrm{and}\;\;\frac{1}{q_{2}}=\frac{1}{p_{2}}+\frac{\eta_{1}}{p'}
\end{align*}
for some
\begin{align*}
-\frac{\tilde{s} p'}{n}<\eta_{1}<\frac{1}{(3-m)(p-1)}.
\end{align*}
 We also take
\begin{align*}
s_{1}=\frac{n}{p_{1}}-\delta_{1}\;\;\textrm{and}\;\;\delta_{3}=\tilde{s}+\frac{n\eta_{1}}{p'}.
\end{align*}
Under present hypotheses, we have $0<s_{1}<(m-1)^{2}n/p_{1}$ and $\delta_{3}>0$. Thus
by the dual principle for $\dot{B}_{p,\infty}^{s_{1}}\times\dot{B}_{p',1}^{-s_{1}}$, inequalities \eqref{est:pro1} and \eqref{ine:mf21},
we can deduce that
\begin{align*}
|I_{2}|&\leq\|J_{m}(u_{1})-J_{m}(u_{2})\|_{\dot{B}_{p_{1},\infty}^{s_{1}}}\|\nabla(v\cdot\varphi)\|_{\dot{B}_{p_{1}',1}^{-s_{1}}}\\
&\leq C\big(\|u_{1}\|_{\dot{B}_{p_{1},1}^{s_{1}/m+(1-1/m)n/p_{1}}}^{m-1}+\|u_{2}\|_{\dot{B}_{p_{1},1}^{s_{1}/m+(1-1/m)n/p_{1}}}^{m-1}\big)\|u_{1}-u_{2}\|_{\dot{B}_{p_{1},1}^{s_{1}/m+(1-1/m)n/p_{1}}}\\
&\quad\cdot\big(\|v\|_{\dot{B}_{p_{2},\infty}^{-\delta_{2}}}\|\varphi\|_{\dot{B}_{p',1}^{1-s_{1}+\delta_{2}}}+\|\varphi\|_{\dot{B}_{q_{1},\infty}^{-\delta_{3}}}\|v\|_{\dot{B}_{q_{2},1}^{1-s_{1}+\delta_{3}}}\big)\\
&\leq C\big(\|u_{1}\|_{\dot{B}_{p,1}^{s_{1}/m+n/p-n/mp_{1}}}^{m-1}+\|u_{2}\|_{\dot{B}_{p,1}^{s_{1}/m+n/p-n/mp_{1}}}^{m-1}\big)\|u_{1}-u_{2}\|_{\dot{B}_{p,1}^{s_{1}/m+n/p-n/mp_{1}}}\\
&\quad\cdot\big(\|v\|_{\dot{B}_{p,\infty}^{n/p_{1}-\delta_{2}}}\|\varphi\|_{\dot{B}_{p',1}^{1-s_{1}+\delta_{2}}}+\|\varphi\|_{\dot{B}_{p',\infty}^{\eta_{1} n/p'-\delta_{3}}}\|v\|_{\dot{B}_{p,1}^{1-s_{1}+\delta_{3}+n/p_{1}-n\eta_{1}/p'}}\big)\\
&\leq C(\|u_{1}\|_{\dot{B}_{p,1}^{s+2\alpha}}^{m-1}+\|u_{2}\|_{\dot{B}_{p,1}^{s+2\alpha}}^{m-1})\|u_{1}-u_{2}\|_{\dot{B}_{p,1}^{s+2\alpha}}\|v\|_{\dot{B}_{p,1}^{s+2\alpha}}\|\varphi\|_{\dot{B}_{p',1}^{-\tilde{s}}}.
\end{align*}

When $2\leq m<\infty$, we use hypotheses $H_{2}$ to deduce that
\begin{align*}
-\frac{n}{p}<\tilde{s}<0\;\;\textrm{and}\;\;0<\frac{n}{p}-\delta_{1}<\min\big\{m-1,\frac{n}{p}\big\}.
\end{align*}
Then similar to the case $1<m<2$, we have
\begin{align*}
|I_{1}|&\leq C\big(\|J_{m}(u_{1})-J_{m}(u_{2})\|_{\dot{B}_{p,\infty}^{n/p-\delta_{1}}}\| v\|_{\dot{B}_{p,\infty}^{\tilde{s}+1+\delta_{1}}}\\
&\quad+\|J_{m}(u_{1})-J_{m}(u_{2})\|_{\dot{B}_{p,\infty}^{\tilde{s}+1+\delta_{1}/m}}\|v\|_{\dot{B}_{p,\infty}^{n/p-\delta_{1}/m}}\big)\|\varphi\|_{\dot{B}_{p',1}^{-\tilde{s} }}\\
&\leq C\big[\big(\|u_{1}\|_{\dot{B}_{p,1}^{n/p-\delta_{1}/m}}^{m-1}+\|u_{2}\|_{\dot{B}_{p,1}^{n/p-\delta_{1}/m}}^{m-1}\big)\|u_{1}-u_{2}\|_{\dot{B}_{p,1}^{n/p-\delta_{1}/m}}\|v\|_{\dot{B}_{p,1}^{\tilde{s}+1+\delta_{1}}}\\
&\quad+\big(\|u_{1}\|_{\dot{B}_{p,1}^{s+2\alpha}}^{m-1}+\|u_{2}\|_{\dot{B}_{p,1}^{s+2\alpha}}^{m-1}\big)\|u_{1}-u_{2}\|_{\dot{B}_{p,1}^{s+2\alpha}}\|v\|_{\dot{B}_{p,\infty}^{n/p-\delta_{1}/m}}\big]\|\varphi\|_{\dot{B}_{p',1}^{-\tilde{s}}}\\
&\leq C(\|u_{1}\|_{\dot{B}_{p,1}^{s+2\alpha}}^{m-1}+\|u_{2}\|_{\dot{B}_{p,1}^{s+2\alpha}}^{m-1})\|u_{1}-u_{2}\|_{\dot{B}_{p,1}^{s+2\alpha}}\|v\|_{\dot{B}_{p,1}^{s+2\alpha}}\|\varphi\|_{\dot{B}_{p',1}^{-\tilde{s} }}.
\end{align*}
where
\begin{align*}
\frac{n}{p}-\delta_{1}=\tilde{s}+1+\frac{\delta_{1}}{m}\;\;\textrm{and}\;\;s+2\alpha=\frac{1}{m}\big(\tilde{s}+1+\frac{\delta_{1}}{m}\big)+\big(1-\frac{1}{m}\big)\frac{n}{p}
\end{align*}
are both taken into account.

Additionally, by taking
\begin{align*}
p_{3}=\frac{p}{p-2}\;\;\textrm{and}\;\;\delta_{4}=\tilde{s}+\frac{n}{p},
\end{align*}
we can also deduce that
\begin{align*}
|I_{2}|&\leq\|J_{m}(u_{1})-J_{m}(u_{2})\|_{\dot{B}_{p,\infty}^{n/p-\delta_{1}}}\|\nabla(v\cdot\varphi)\|_{\dot{B}_{p',1}^{\delta_{1}-n/p}}\\
&\leq C\big(\|u_{1}\|_{\dot{B}_{p,1}^{n/p-\delta_{1}/m}}^{m-1}+\|u_{2}\|_{\dot{B}_{p,1}^{n/p-\delta_{1}/m}}^{m-1}\big)\|u_{1}-u_{2}\|_{\dot{B}_{p,1}^{n/p-\delta_{1}/m}}\\
&\quad\cdot\big(\|v\|_{\dot{B}_{\infty,\infty}^{-\delta_{1}/m}}\|\varphi\|_{\dot{B}_{p',1}^{1-s_{2}+\delta_{1}/m}}+\|\varphi\|_{\dot{B}_{p_{3},\infty}^{-\delta_{4}}}\|v\|_{\dot{B}_{p,1}^{1+\delta_{1}-n/p+\delta_{4}}}\big)\\
&\leq C\big(\|u_{1}\|_{\dot{B}_{p,1}^{n/p-\delta_{1}/m}}^{m-1}+\|u_{2}\|_{\dot{B}_{p,1}^{n/p-\delta_{1}/m}}^{m-1}\big)\|u_{1}-u_{2}\|_{\dot{B}_{p,1}^{n/p-\delta_{1}/m}}\\
&\quad\cdot\big(\|v\|_{\dot{B}_{p,\infty}^{n/p-\delta_{1}/m}}\|\varphi\|_{\dot{B}_{p',1}^{1-s_{2}+\delta_{1}/m}}+\|\varphi\|_{\dot{B}_{p',\infty}^{n/p-\delta_{4}}}\|v\|_{\dot{B}_{p,1}^{1+\delta_{1}-n/p+\delta_{4}}}\big)\\
&\leq C(\|u_{1}\|_{\dot{B}_{p,1}^{s+2\alpha}}^{m-1}+\|u_{2}\|_{\dot{B}_{p,1}^{s+2\alpha}}^{m-1})\|u_{1}-u_{2}\|_{\dot{B}_{p,1}^{s+2\alpha}}\|v\|_{\dot{B}_{p,1}^{s+2\alpha}}\|\varphi\|_{\dot{B}_{p',1}^{-\tilde{s}}},
\end{align*}
where $\delta_{4}>0$ is assured by present hypotheses.

Based on the above derivations, with the aid of the density of $\mathcal{S}_{0}$ in both $\dot{B}_{p,1}^{s+2\alpha}$ and $\dot{B}_{p',1}^{-\tilde{s}}$, we conclude that
\begin{align*}
&\quad|\langle J_{m}(u_{1}(t))-J_{m}(u_{2}(t)))\cdot\nabla v(t),\varphi\rangle|\\&\leq C(\|u_{1}(t)\|_{\dot{B}_{p,1}^{s+2\alpha}}^{m-1}+\|u_{2}(t)\|_{\dot{B}_{p,1}^{s+2\alpha}}^{m-1})\|u_{1}(t)-u_{2}(t)\|_{\dot{B}_{p,1}^{s+2\alpha}}\|v(t)\|_{\dot{B}_{p,1}^{s+2\alpha}}\|\varphi\|_{\dot{B}_{p',1}^{-\tilde{s}}},
\end{align*}
for all $\varphi\in\dot{B}_{p',1}^{-\tilde{s}}$ and consequently
\begin{align*}
&\quad\|J_{m}(u_{1}(t))-J_{m}(u_{2}(t)))\cdot\nabla v(t)\|_{\dot{B}_{p,\infty}^{\tilde{s}}}\\&\leq C(\|u_{1}(t)\|_{\dot{B}_{p,1}^{s+2\alpha}}^{m-1}+\|u_{2}(t)\|_{\dot{B}_{p,1}^{s+2\alpha}}^{m-1})\|u_{1}(t)-u_{2}(t)\|_{\dot{B}_{p,1}^{s+2\alpha}}\|v(t)\|_{\dot{B}_{p,1}^{s+2\alpha}},
\end{align*}
for a.e. $t\in(0,T)$.

Finally by invoking Lemma \ref{lem:pfg}, \ref{lem:pm}, one can derive that
\begin{align*}
&\quad\|(J_{m}(u_{1})-J_{m}(u_{2}))\cdot\nabla v\|_{L^{\tilde{\rho},r}(0,T;\dot{B}_{p,\infty}^{\tilde{s}})}\\
&\leq C\big\|\|u_{1}\|_{\dot{B}_{p,1}^{s+2\alpha}}^{m-1}+\|u_{2}\|_{\dot{B}_{p,1}^{s+2\alpha}}^{m-1}\big\|_{L^{\rho/(m-1),r}(0,T)}\\
&\quad\cdot\big\|\|u_{1}-u_{2}\|_{\dot{B}_{p,1}^{s+2\alpha}}\big\|_{L^{\rho,r}(0,T)}\big\|\|v\|_{\dot{B}_{p,1}^{s+2\alpha})}\big\|_{L^{\rho ,r}(0,T)}\\
&\leq C(\|u_{1}\|_{L_{T}^{\rho,r}(\dot{B}_{p,1}^{s+2\alpha})}^{m-1}+\|u_{2}\|_{L_{T}^{\rho,r}(\dot{B}_{p,1}^{s+2\alpha})}^{m-1})\|u_{1}-u_{2}\|_{L_{T}^{\rho,r}(\dot{B}_{p,1}^{s+2\alpha})}\|v\|_{L_{T}^{\rho,r}(\dot{B}_{p,1}^{s+2\alpha})}.
\end{align*}
Thus the proof has been completed. \hfill $\Box$

\noindent{\it Proof of Theorem \ref{thm:1.1}} :
Introduce a nonlinear operator
\begin{equation}\label{eqn:phi}
  \Phi(u):=a_{L}+\mathcal{R}(\mathcal{P}f-\mathcal{P}(J_{m}(u)\cdot\nabla u))
\end{equation}
Evidently by \eqref{rep:u}, each strong solution of \eqref{eqn:gns1} is a fixed point of $\Phi$ in some suitable function space.

We now invoking inequalities \eqref{ine:exp1}, \eqref{est:f1} and Proposition \ref{prop:3.1}, \ref{prop:3.2}, together with the boundedness of $\mathcal{P}$ on Besov spaces to make estimates for $a_{L}$ and $\mathcal{R}(\mathcal{P}f-\mathcal{P}(J_{m}(u)\cdot\nabla u))$ under hypotheses $H_{i}$, $i=0,1,2$ for $m=1$, $1<m<2$ or $2\leq m<\infty$ respectively.
\begin{align}\label{est:ul}
\|a_{L}\|_{L_{T}^{\rho,r}(\dot{B}_{p,1}^{s+2\alpha})}\leq k_{0}\|a\|_{\dot{B}_{p_{0},r}^{s_{0}}},
\end{align}
\begin{align}\label{est:sfu1}
\nonumber&\quad\|\mathcal{R}(\mathcal{P}f-\mathcal{P}(J_{m}(u)\cdot\nabla u))\|_{L_{T}^{\rho,r}(\dot{B}_{p,1}^{s+2\alpha})}\\
&\leq k_{1}\big(\|f\|_{L_{T}^{\tilde{\rho},r}(\dot{B}_{p,\infty}^{\tilde{s}})}+\|J_{m}(u)\cdot\nabla u\|_{L_{T}^{\tilde{\rho},r}(\dot{B}_{p,\infty}^{\tilde{s}})}\big)\\
\nonumber&\leq k_{1}\|f\|_{L_{T}^{\tilde{\rho},r}(\dot{B}_{p,\infty}^{\tilde{s}})}+k_{2}\|u\|_{L_{T}^{\rho,r}(\dot{B}_{p,1}^{s+2\alpha})}^{m+1},
\end{align}
and
\begin{align}\label{est:sfu2}
\nonumber&\quad\|\mathcal{R}(\mathcal{P}(J_{m}(u_{1})\cdot\nabla u_{1}))-\mathcal{R}(\mathcal{P}(J_{m}(u_{2})\cdot\nabla u_{2}))\|_{L_{T}^{\rho,r}(\dot{B}_{p,1}^{s+2\alpha})}\\
\nonumber&\leq k_{1}\big(\|(J_{m}(u_{1})-J_{m}(u_{2}))\cdot\nabla u_{2}\|_{L_{T}^{\tilde{\rho},r}(\dot{B}_{p,\infty}^{\tilde{s}})}\\
&\quad+\|J_{m}(u_{1})\cdot\nabla (u_{1}-u_{2})\|_{L_{T}^{\tilde{\rho},r}(\dot{B}_{p,\infty}^{\tilde{s}})}\big)\\
\nonumber&\leq k_{2}(\|u_{1}\|_{L_{T}^{\rho,r}(\dot{B}_{p,1}^{s+2\alpha})}^{m-1}+\|u_{2}\|_{L_{T}^{\rho,r}(\dot{B}_{p,1}^{s+2\alpha})}^{m-1})\|u_{1}-u_{2}\|_{L_{T}^{\rho,r}(\dot{B}_{p,1}^{s+2\alpha})}\|u_{2}\|_{L_{T}^{\rho,r}(\dot{B}_{p,1}^{s+2\alpha})}\\
\nonumber&\quad+k_{2}\|u_{1}\|_{L_{T}^{\rho,r}(\dot{B}_{p,1}^{s+2\alpha})}^{m}\|u_{1}-u_{2}\|_{L_{T}^{\rho,r}(\dot{B}_{p,1}^{s+2\alpha})},
\end{align}
where the constants $k_{i}\geq1$, $i=0,1,2$ independent of $a,f$ and $u$, $u_{i}$, $i=1,2$ are fixed for later use.

Plugging \eqref{est:ul}--\eqref{est:sfu2} into \eqref{eqn:phi}, we obtain
\begin{align}\label{est:phu1}
\nonumber\|\Phi(u)\|_{L_{T}^{\rho,r}(\dot{B}_{p,1}^{s+2\alpha})}&\leq k_{0}\|a\|_{\dot{B}_{p_{0},r}^{s_{0}}}+k_{1}\|f\|_{L_{T}^{\tilde{\rho},r}(\dot{B}_{p,\infty}^{\tilde{s}})}+k_{2}\|u\|_{L_{T}^{\rho,r}(\dot{B}_{p,1}^{s+2\alpha})}^{m+1}\\
&=K_{0}+k_{2}\|u\|_{L_{T}^{\rho,r}(\dot{B}_{p,1}^{s+2\alpha})}^{m+1}
\end{align}
and
\begin{align}\label{est:phu2}
\nonumber&\quad\|\Phi(u_{1})-\Phi(u_{2})\|_{L_{T}^{\rho,r}(\dot{B}_{p,1}^{s+2\alpha})}\\
&\leq 2k_{2}(\|u_{1}\|_{L_{T}^{\rho,r}(\dot{B}_{p,1}^{s+2\alpha})}^{m}+\|u_{2}\|_{L_{T}^{\rho,r}(\dot{B}_{p,1}^{s+2\alpha})}^{m})\|u_{1}-u_{2}\|_{L_{T}^{\rho,r}(\dot{B}_{p,1}^{s+2\alpha})},
\end{align}
where the constant
\begin{align*}
K_{0}:=k_{0}\|a\|_{\dot{B}_{p_{0},r}^{s_{0}}}+k_{1}\|f\|_{L_{T}^{\tilde{\rho},r}(\dot{B}_{p,\infty}^{\tilde{s}})}
\end{align*}
is assumed to less than $\eta:=1/16k_{2}$ so that the quadratic equation $k_{2}\lambda^{2}-\lambda+K_{0}=0$ has two different positive root, the first of which, named $\lambda_{1}=(1-\sqrt{1-4k_{2}K_{0}})/2k_{2}$, is less than $1/4k_{2}$.

Denote by $B(0,\lambda_{1})$ the closed ball in $L_{T}^{\rho,r}(\dot{B}_{p,1}^{s+2\alpha})$ centered at 0 with the radius $\lambda_{1}$. From \eqref{est:phu1} and the fact $\lambda_{1}<1$ and $m>1$, it is easy to see that if $\|u\|_{L_{T}^{\rho,r}(\dot{B}_{p,1}^{s+2\alpha})}\leq\lambda_{1}$, then
\begin{align*}
 \|\Phi(u)\|_{L_{T}^{\rho,r}(\dot{B}_{p,1}^{s+2\alpha})}\leq K_{0}+k_{2}\lambda_{1}^{m+1}\leq \lambda_{1}.
\end{align*}
Moreover, if $\|u_{i}\|_{L_{T}^{\rho,r}(\dot{B}_{p,1}^{s+2\alpha})}\leq\lambda_{1}$, $i=1,2$, then
\begin{equation}\label{ine:u12}
  \|\Phi(u_{1})-\Phi(u_{2})\|_{L_{T}^{\rho,r}(\dot{B}_{p,1}^{s+2\alpha})}\leq 4k_{2}\lambda_{1}\|u_{1}-u_{2}\|_{L_{T}^{\rho,r}(\dot{B}_{p,1}^{s+2\alpha})}.
\end{equation}
Since $4k_{2}\lambda_{1}<1$, we conclude that $\Phi:B(0,\lambda_{1})\rightarrow B(0,\lambda_{1})$ is a contraction. Consequently, it has a unique fixed point $u\in B(0,\lambda_{1})$, which verifies the integral equation
\begin{equation}\label{eq:u1}
  u=a_{L}+\mathcal{R}(\mathcal{P}f-\mathcal{P}(J_{m}(u)\cdot\nabla u)),
\end{equation}
and obeys the following estimate
\begin{equation}\label{est:slu}
  \|u\|_{L_{T}^{\rho,r}(\dot{B}_{p,\infty}^{s+2\alpha})}\leq2K_{0}=2\big(k_{0}\|a\|_{\dot{B}_{p_{0},r}^{s_{0}}}+k_{1}\|f\|_{L_{T}^{\tilde{\rho},r}(\dot{B}_{p,\infty}^{\tilde{s}})}\big).
\end{equation}

We next show that $u$ is the strong solution of \eqref{eqn:gns1}. For this purpose, we notice the following facts: Both $\mathcal{P}f$ and $\mathcal{P}(J_{m}(u)\cdot\nabla u))$ belong to $L_{T}^{\tilde{\rho},r}(\dot{B}_{p,\infty}^{\tilde{s}})$, the initial velocity $a$ lies in $\dot{B}_{p_{0},r}^{s_{0}}$, where the associated exponents verify
\begin{align*}
s_{0}-\frac{n}{p_{0}}-2\alpha=s-\frac{n}{p}-\frac{2\alpha}{\rho}=\tilde{s}-\frac{n}{p}-\frac{2\alpha}{\tilde{\rho}}.
\end{align*}
Thus in light of Proposition \ref{prop:1.1a}, we assert that initial problem
 \begin{align*}
\partial_{t}v+\mu(-\Delta)^{\alpha}v=-\mathcal{P}(J_{m}(u)\cdot\nabla u)+\mathcal{P}f(t),\;\;t>0, \;\; u(0)=a.
 \end{align*}
has a unique strong solution, say $v$, lying in the class $L_{T}^{\tilde{\rho},r}(\dot{B}_{p,\infty}^{\tilde{s}+2\alpha})$ with $\partial_{t}v\in L_{T}^{\tilde{\rho},r}(\dot{B}_{p,\infty}^{\tilde{s}})$. This solution is exactly equal to $u$, since both of them solve the integral equation \eqref{eq:u1} with the right side identical. Moreover, it is easy to check that $a_{L}\in C([0,T];\dot{B}_{p_{0},r}^{s_{0}})$ and $\mathcal{R}(\mathcal{P}f-\mathcal{P}(J_{m}(u)\cdot\nabla u))\in C([0,T];\dot{B}_{p,\infty}^{\tilde{s}})$, so as their sum, function $u$ lies in $C([0,T];\dot{B}_{p_{0},r}^{s_{0}}+\dot{B}_{p,\infty}^{\tilde{s}})$ together with $u(0)=a$. Conclusion, $u$ is the strong solution of \eqref{eqn:gns1}.

By invoking Proposition \ref{prop:1.1a} and inequalities \eqref{est:uv1}, \eqref{est:slu}, we further get estimates for higher regularity of $u$, that is
\begin{align*}
&\quad\|u\|_{L_{T}^{\tilde{\rho},r}(\dot{B}_{p,\infty}^{\tilde{s}+2\alpha})}+\|\partial_{t}u\|_{L_{T}^{\tilde{\rho},r}(\dot{B}_{p,\infty}^{\tilde{s}})}\\
&\leq C\big(\|a\|_{\dot{B}_{p_{0},r}^{s_{0}}}+\|\mathcal{P}f\|_{L_{T}^{\tilde{\rho},r}(\dot{B}_{p,\infty}^{\tilde{s}})}+\|\mathcal{P}(J_{m}(u)\cdot\nabla u)\|_{L_{T}^{\tilde{\rho},r}(\dot{B}_{p,\infty}^{\tilde{s}})} \big)\\
&\leq C\big(\|a\|_{\dot{B}_{p_{0},r}^{s_{0}}}+\|f\|_{L_{T}^{\tilde{\rho},r}(\dot{B}_{p,\infty}^{\tilde{s}})}+\|u\|_{L_{T}^{\rho,r}(\dot{B}_{p,\infty}^{s+2\alpha})}^{m+1} \big)\\
&\leq C\big(\|a\|_{\dot{B}_{p_{0},r}^{s_{0}}}+\|f\|_{L_{T}^{\tilde{\rho},r}(\dot{B}_{p,\infty}^{\tilde{s}})}+\|a\|_{\dot{B}_{p_{0},r}^{s_{0}}}^{m+1}+\|f\|_{L_{T}^{\tilde{\rho},r}(\dot{B}_{p,\infty}^{\tilde{s}})}^{m+1}\big)\\
&\leq C\big(\|a\|_{\dot{B}_{p_{0},r}^{s_{0}}}+\|f\|_{L_{T}^{\tilde{\rho},r}(\dot{B}_{p,\infty}^{\tilde{s}})}\big).
 \end{align*}
In addition, if we let
\begin{align*}
\nabla \pi=(I-\mathcal{P})(f-J_{m}(u)\cdot\nabla u),
\end{align*}
then we have $\nabla \pi\in L_{T}^{\tilde{\rho},r}(\dot{B}_{p,\infty}^{\tilde{s}})$, and
\begin{align*}
\|\nabla \pi\|_{L_{T}^{\tilde{\rho},r}(\dot{B}_{p,\infty}^{\tilde{s}})}\leq C\big(\|a\|_{\dot{B}_{p_{0},r}^{s_{0}}}+\|f\|_{L_{T}^{\tilde{\rho},r}(\dot{B}_{p,\infty}^{\tilde{s}})}\big).
\end{align*}
Thus inequality \eqref{est:sl} has been reached.

It is easy to check that the function couple $(u,\nabla \pi)$ verify the generalised Navier-Stokes equation \eqref{eqn:gns} in the space $\dot{B}_{p,\infty}^{\tilde{s}}$ a.e. on $[0,T]$. So they comprise the strong solution of \eqref{eqn:gns}. Finally, uniqueness of the strong solution comes from inequality \eqref{ine:u12}, assumption \eqref{cd:af} and the formation of $\nabla\pi$. Thus Theorem \ref{thm:1.1} has been proved. \hfill $\Box$

\begin{rem}
Apart from \eqref{hy:p1}, if we assume further $p>m(2-m)n/2(3-m)$, it follows $m(2-m)n/2(3-m)p<1/2$, which means that exponent $\alpha$ fulfilling \eqref{hy:a1} can take the value $1$ in the case $1<m<2$. This fact also holds for other cases. Therefore results of Theorem \ref{thm:1.1} remain valid for the (GNS) \eqref{eqn:gns} with $\alpha=1$ naturally.
\end{rem}


\begin{thebibliography}{99}

\bibitem{ho1951} E. Hopf, \"{U}ber die Anfangswertaufgabe f\"{u}r die hydrodynamischen Grundgleichungen,  Math. Nachr. 4 (1951) 213--231.

\bibitem{ms1988} T. Miyakawa, H. Sohr, On energy inequality, smoothness and large time behavior in $L^{2}$ for weak solutions of the Navier-
Stokes equations in exterior domains,  Math. Z. 199 (1988) 455--478.

\bibitem{ga2000} Galdi, G. P.: An Introduction to the Navier-Stokes Initial-Boundary Value Problem, Fundamental Directions in Mathematical Fluid Mechanics, Birkh\"{a}user Basel, 2000.

\bibitem{fk1964} H. Fujita, T. Kato, On the Navier-Stokes initial problem I, Arch. Rat. Mech
Anal. 16 (1964) 269--315.

 \bibitem{ka1984} T. Kato, Strong $L^{p}-$solutions of the Navier-Stokes equation in $\mathbb{R}^{m}$, with applications to weak solutions,  Math. Z. 187 (1984) 471--480.

 \bibitem{ko1989} H. Kozono, Global $L^{n}-$solution and its decay property for the Navier-Stokes equations in half-space $\mathbb{R}_{+}^{n}$,  J. Differential  Equations 79 (1989) 79--88.

 \bibitem{iw1989} Iwashita, H.: $L_{q}-L_{r}$ estimates for solutions of the nonstationary Stokes equations in an exterior domain and the Navier-Stokes initial value Problems in $L_{q}$ spaces, Math. Ann. 285 (1989) 265--288.

\bibitem{li1969} J. L. Lions, Quelques M\'{e}thodes de R\'{e}solution des Probl\'{e}mes aux Limites Non Lin\'{e}aires, Vol 1. Dunod, Paris, 1969.

\bibitem{ch2017}  Z-M. Chen, Analitic semigroup approach to generalized Navier-Stokes flows in Besov spaces, J. Math.
Fluid Mech. 19 (2017) 709--724.

 \bibitem{du2018}  N. Duan, Well-posedness and decay of solutions for three dimensional  generalized Navier-Stokes equations, Comput.  Math. Appl. 76 (2018) 1026--1033.

\bibitem{li2020} Q. Liu, Global well-posedness of the generalized incompressible
Navier-Stokes equations with large initial data, Bull. Malays. Math. Sci. Soc. 43(2) (2020) 2549--2564.

\bibitem{wu2005} J. Wu, Lower bounds for an integral involving fractional Laplacians and the generalized Navier-Stokes equations in Besov
spaces, Comm. Math. Phys. 263 (2005) 803--831.

\bibitem{yz2014} X. Yu, Z. Zhai, Well-posedness for fractional Navier-Stokes equations in the largest critical spaces $\dot{B}_{\infty,\infty}^{-(2\beta-1)}(\mathbb{R}^{n})$, Math. Methods Appl. Sci. 35 (2014) 676--683.

\bibitem{ks2019} H. Kozono, S. Shimizu, Strong solutions of the Navier-Stokes equations based on the maximal Lorentz regularity theorem in Besov spaces, J. Func. Anal. 276 (2019) 896--931.

\bibitem{bcd2011} H. Bahouri, J.-Y. Chemin, R. Danchin, Fourier analysis and Nonlinear Partial Differential Equations.
A Series of Comprehensive Studies in Mathematics, vol 343,  Springer-Verlag, 2011.

\bibitem{mwz2012} C. Miao, J. Wu, Z. Zhang, Littlewood Paley Decomposition Theory and Its Application in Fluid Dynamics Equations (in Chinese), Foundamental Series of Modern Mathematics, vol 142, Science Press, 2012.

\bibitem{zcz2024} Q. Zhang, S. Cao \& H. Zhang, Global solvability of the generalised Boussinesq system with linear or nonlinear buoyancy force, Submitted for pubication.

\end{thebibliography}
\end{document}